\documentclass[11pt,a4paper]{article}
\usepackage{latexsym,amssymb,amsthm,amsmath}
\title{On the structure of Hardy-Sobolev-Maz'ya inequalities}

\author{\Large Stathis Filippas$^{1,4}$,  Achilles  Tertikas$^{2,4}$
\& Jesper Tidblom$^{3}$  \\
                                                                           \\
        Department of Applied Mathematics$^{1}$ \\
         University of Crete,
         71409 Heraklion,  Greece \\
        filippas@tem.uoc.gr\\
                  \\
  Department of Mathematics$^{2}$ \\
     University of Crete,
         71409 Heraklion,  Greece          \\
  tertikas@math.uoc.gr\\
\\
    The Erwin Schr\"odinger Institute (ESI)$^{3}$ \\
          Boltzmanngasse 9, A-1090 Vienna, Austria \\
          Jesper.Tidblom@esi.ac.at\\
                                      \\
        Institute of Applied and Computational Mathematics$^4$, \\
        FORTH, 71110 Heraklion, Greece \\
    \\ }

\date{\today}
\newtheorem{theorem}{Theorem}[section]
\newtheorem{proof*}{Proof:}
\newtheorem{corollary}[theorem]{Corollary}
\newtheorem{lemma}[theorem]{Lemma}
\numberwithin{equation}{section}
\newcommand{\be}{\begin{equation}}
\newcommand{\ee}{\end{equation}}
\newcommand{\bea}{\begin{eqnarray}}
\newcommand{\eea}{\end{eqnarray}}
\newcommand{\la}{\label}
\newcommand{\xa}{\alpha}
\newcommand{\xb}{\beta}
\newcommand{\xg}{\gamma}
\newcommand{\xe}{\varepsilon}
\newcommand{\xs}{\sigma}

\newcommand{\ra}{\rightarrow}
\def\finedim{$\hfill \Box$}

\begin{document}
\maketitle

\begin{abstract}
 \noindent
In this article we  establish new improvements of the optimal
Hardy inequality in the half space. We first
 add all possible linear combinations of
Hardy type terms thus revealing the structure of this type
of inequalities  and obtaining best constants. We then add
 the  critical Sobolev term  and obtain necessary and sufficient
conditions for the validity of  Hardy-Sobolev-Maz'ya type
inequalities.
\end{abstract}

\section{Introduction}
One version of the Hardy inequality   states that for convex
domains  $\Omega\subset\mathbb{R}^n$  the following estimate holds
\[
 \int_{\Omega}|\nabla u|^2dx \geq  \frac{1}{4}
 \int_{\Omega}\frac{|u|^2}{d(x)^2}dx, ~~~~~~~~~~~ u \in
 C^{\infty}_0(\Omega),
\]
where $d(x)=\mathrm{dist}(x,\partial \Omega)$ and the constant
$\frac{1}{4}$ is the best possible constant. This result has been
improved and generalized in many different ways, see for example
\cite{Ana},  \cite{Bar2}, \cite{Bre}, \cite{Dav}, \cite{DD},
\cite{Fil}, \cite{Fil2},
 \cite{Hoff2}, \cite{Tid1}, \cite{Tid2}.

One pioneering result due to Brezis and Marcus \cite{Bre}
  is the following improved
Hardy inequality:
\begin{equation}\label{1.1}
\int_{\Omega}|\nabla u|^2dx \geq  \frac{1}{4}
 \int_{\Omega}\frac{|u|^2}{d(x)^2}dx
 + C_2(\Omega) \int_{\Omega}|u|^2 dx, \quad u
 \in C^{\infty}_0(\Omega),
\end{equation}
valid for  any convex domain $\Omega\subset\mathbb{R}^n$. This
estimate has been recently extended in \cite{Fil2}:
\begin{equation}
\int_{\Omega}|\nabla u|^2dx \geq  \frac{1}{4}
 \int_{\Omega}\frac{|u|^2}{d(x)^2}dx
 + C_q(\Omega)\left( \int_{\Omega}|u|^qdx \right)^{\frac{2}{q}}, \quad u
 \in C^{\infty}_0(\Omega). \label{geometric}
\end{equation}
Moreover, it is shown in \cite{Fil2} that
 there
exist constants $c_1$ and $c_2$ only depending on $q$ and the
dimension $n$ of $\Omega$ such that the best constant $ C_q(\Omega)$
satisfies
\[
 c_1 D^{n-2-\frac{2n}{q}} \geq C_q(\Omega) \geq
 c_2D^{n-2-\frac{2n}{q}},
\]
where $D= \mathrm{sup}_{x\in\Omega}d(x) < \infty$ and $2 \leq q <
\frac{2n}{n-2}$.
 We note that the critical  Sobolev  exponent  $q=2^{*}:=\frac{2n}{n-2}$
is not included in the above theorem. For results in the critical
case  we refer to \cite{Fil}.

  Let us denote by $S_n=
\pi n(n-2) \left( \frac{\Gamma(\frac{n}{2})}{\Gamma(n)}
 \right)^{2/n}$, $n \geq 3$,
the best constant in the Sobolev inequality
\[
   \int_{\Omega}|\nabla u|^2dx  \geq S_n
 \left( \int_{\Omega}|u|^{2^*}dx \right)^{\frac{2}{2^*}},
   \quad u \in C_0^{\infty}(\Omega).
\]
The first inequality that combines both the critical Sobolev exponent
term and the Hardy  term the latter with best constant, is due to
Maz'ya \cite{Maz}, and is the following Hardy--Sobolev--Maz'ya
inequality:
\begin{equation}
\int_{\mathbb{R}^n_+}|\nabla u|^2dx \geq \frac{1}{4}
\int_{\mathbb{R}^n_+}\frac{|u|^2}{x_1^2}dx + C_n \left(
\int_{\mathbb{R}^n_+}|u|^{2^*}dx \right)^{\frac{2}{2^*}}, \quad u
\in C^{\infty}_0(\mathbb{R}^n_+), \label{maz}
\end{equation}
where $\mathbb{R}^n_+ = \{ (x_1, \ldots, x_n): x_1 >0 \}$
 denotes the upper half-space, $C_n$ is a positive
 constant and $2^*=2n/(n-2)$,  $n \geq 3$. Recently, it was shown in \cite{BFL}
 that
in the 3--dimensional case $n=3$, the best constant $C_3$ coincides
with the best Sobolev constant  $S_3$! On the other hand when
$n \geq 4$ one  has that $C_n < S_n$,  see \cite{TT}.

We next mention an improvement of  Hardy's inequality
 that involves two distance
functions:
\[
\int_{\mathbb{R}^n_+}|\nabla u|^2dx \geq  \frac14
 \int_{\mathbb{R}^n_+}\frac{|u|^2}{x_1^2}dx + C(\tau)
\int_{\mathbb{R}^n_+}\frac{|u|^2}{x_1^{2-\tau}
(x_{1}^2+x_{2}^2)^{\frac{\tau}{2}}}dx,
\qquad u
\in C^{\infty}_0(\mathbb{R}^n_+),
\]
where $0<\tau  \leq 1$.
This is a special case of a  more general inequality proved in
 \cite{Tid2}.

In this work we study   improvements of Hardy's inequality
that involve various distance functions.
Working in the upper half space $\mathbb{R}^n_+$,
we obtain Hardy type inequalities that involve constant multiples
of the inverse square of  the distance
to
linear submanifolds of different codimensions of the boundary
 $\partial \mathbb{R}^n_+$.
Actually, we are able to give a complete description of the
structure of this kind of improved Hardy  inequalities.
 In particular,
 we
have a lot of freedom in
choosing these constants and we will show that all our configurations of constants are, in a natural sense, optimal.
More precisely, our first result reads:

\quad
\\ \textbf{Theorem A}
     ({\bf Improved Hardy inequality}) \\
 {\it {\bf i)} Let $\alpha_1, \alpha_2,
\ldots, \alpha_n$  be arbitrary real numbers and
\begin{eqnarray*}
 \beta_1 &=& -\alpha_1^2   +\frac14,  \\
 \beta_m  &=&  - \alpha_m^2 + \left(\alpha_{m-1}-\frac12 \right)^2,
~~~~m=2,3,\ldots,n.
\end{eqnarray*}
Then for any $u \in C_0^{\infty}(\mathbb{R}^n_+)$ there holds
\[
  \int_{\mathbb{R}^n_+}|\nabla u|^2dx \geq \int_{\mathbb{R}^n_+}\left(\frac{\beta_1}{x_1^2}+
 \frac{\beta_2}{x_1^2+x_2^2} + \ldots + \frac{\beta_n}{x_1^2+x_2^2+\ldots+x_n^2} \right)u^2dx.
\]
{\bf ii)}Suppose that for some real numbers $\beta_1, \beta_2
\ldots, \beta_n$ the
 following
inequality holds \be\nonumber
  \int_{\mathbb{R}^n_+}|\nabla u|^2dx \geq
\int_{\mathbb{R}^n_+}\left(\frac{\beta_1}{x_1^2}+
 \frac{\beta_2}{x_1^2+x_2^2} + \ldots
 + \frac{\beta_n}{x_1^2+x_2^2+\ldots+x_n^2} \right)u^2dx,
\ee for any $u \in C_0^{\infty}(\mathbb{R}^n_+)$. Then, there
exists  nonpositive constants
 $\xa_1, \ldots, \xa_n$, such that
\begin{eqnarray*}
 \beta_1 &=& -\alpha_1^2   +\frac14,  \\
 \beta_m  &=&  - \alpha_m^2 + \left(\alpha_{m-1}-\frac12 \right)^2,
~~~~m=2,3,\ldots,n.
\end{eqnarray*}
}

We next investigate  the possibility of adding Sobolev
type remainder terms.  It
turns out that almost every choice of the constants in theorem A
allows one to add a positive Sobolev term as well. The details are
in our second main theorem.

\quad \\ \textbf{Theorem B}
({\bf Improved Hardy--Sobolev--Maz'ya  inequality}) \\
 {\it Let $\alpha_1, \alpha_2, \ldots,
\alpha_n$  be arbitrary nonpositive
 real numbers and
\begin{eqnarray*}
 \beta_1 &=& -\alpha_1^2   +\frac14,  \\
 \beta_m  &=&  - \alpha_m^2 + \left(\alpha_{m-1}-\frac12 \right)^2,
~~~~m=2,3,\ldots,n.
\end{eqnarray*}
Then,  if $\xa_n <0$ there exists a positive constant $C$ such
that
 for any
$u \in C_0^{\infty}(\mathbb{R}^n_+)$ there holds \be\la{1.20}
  \int_{\mathbb{R}^n_+}|\nabla u|^2dx \geq \int_{\mathbb{R}^n_+}\left(\frac{\beta_1}{x_1^2}+
 \frac{\beta_2}{x_1^2+x_2^2} + \ldots + \frac{\beta_n}{x_1^2+x_2^2+\ldots+x_n^2} \right)u^2dx + C\left(\int_{\mathbb{R}^n_+}|u|^{2^*}dx
  \right)^{\frac{2}{2^*}}.
\ee If  $\xa_n=0$ then there is no positive constant $C$ such that
(\ref{1.20}) holds. }
\\ \quad
\\
It is interesting to note that the Sobolev term vanishes precisely
when the constant $\beta_n$, in front of the Hardy-type term
containing the point singularity, is chosen optimal. It is a bit
curious that the size of the other constants,
$\beta_1,\ldots,\beta_{n-1}$, does not matter at all for this
question. Only the relative size of $\beta_n$ compared to the
other constants matters.

Our results depend heavily on the Gagliardo-Nirenberg-Sobolev
inequality and also on an interesting relation between the
existence of an $L^1$ Hardy inequality and the possibility of  adding
a Sobolev type remainder term to the corresponding $L^2$ inequality.
The precise result reads:
 \\
\quad \\ \textbf{Theorem C} {\it
Let $\Omega \subset  \mathbb{R}^n$, $n \geq 3$ be  a smooth domain.
 Assume that $\phi > 0$, $\phi \in
C^2(\Omega)$ and that the following weighted $L^1$
inequality holds
\begin{equation}\la{c1}
 \int_{\Omega} \phi^{\frac{2(n-1)}{n-2}}|\nabla
 v|dx \geq C\int_{\Omega} \phi^{\frac{n}{n-2}}|\nabla \phi||v|dx ,
 \quad v \in C^{\infty}_0(\Omega).
\end{equation}
Then, there exists $c>0$ such  that
\be\la{c2}
 \int_{\Omega}|\nabla u|^2 dx
 \geq -\int_{\Omega}\frac{\Delta \phi}{\phi}|u|^2dx + c
 \left( \int_{\Omega} |u|^{2^*}dx
 \right)^{\frac{2}{2^*}}, ~~~~ u \in C^{\infty}_0(\Omega).
\ee}
 The regularity  assumptions  on $\phi$  can be weakened,
 but for our purposes it is enough to restrict
ourselves to $\phi \in C^2(\Omega)$. We note  that under
the sole assumption $\phi>0$ and $\phi \in C^2(\Omega)$  the
following inequality
\be\label{in20}
\int_{\Omega}|\nabla u|^2 dx
 \geq -\int_{\Omega}\frac{\Delta \phi}{\phi}|u|^2dx,
 ~~~~ u \in C^{\infty}_0(\Omega).
\ee
is always true;  see Lemma  \ref{mainthm}.  It is the validity
of (\ref{c1}) that makes possible the  addition of  the Sobolev term in
(\ref{in20}). An easy example where   both (\ref{c1}) and (\ref{c2})
fail, is the
 case where $\phi$ is taken to be the first  Dirichlet eigenfunction
of the Laplacian of $\Omega$, for $\Omega$ bounded.

 Our
methods are not restricted to the case $\Omega = \mathbb{R}^n_+$.
In the last section of the paper we give an example of how to
apply the method to get some results for the quarter-space.
Moreover, as one can easily  check  our results remain valid
even for complex valued functions.

The paper is organized as follows. In section 2 we give the proof
of Theorem A. In section 3 we give the proofs of Theorems B and C.
Finally, in the last section we obtain some results for the quarter
space.

{ \bf Acknowledgment} This work was largely done whilst  JT
 was visiting the University of Crete  and FORTH  in
Heraklion, supported by a postdoctoral  fellowship through the
  RTN  European network
Fronts--Singularities, HPRN-CT-2002-00274.   SF and AT acknowledge partial
support by the same program.

\section{Improved Hardy inequalities in the half-space}

The half-space $\mathbb{R}^n_+$
has some nice features that are not present for an arbitrary
convex domain. The fact that the boundary has zero curvature is
very useful when one is trying to prove certain
sorts of inequalities, as we shall see below.

We start with a general  auxiliary  Lemma.

\begin{lemma} \label{mainthm}
(i) Let $\mathbf{F} \in C^1(\Omega)$, then
\begin{equation} \label{eq1}
\int_{\Omega}|\nabla u|^2 dx = \int_{\Omega}\left( \mathrm{div}
\mathbf{F} - |\mathbf{F}|^{2}\right)|u|^2dx+\int_{\Omega} |\nabla
u+\mathbf{F}u|^2dx, ~~~~ \forall u \in C^{\infty}_0(\Omega).
\end{equation}
(ii) Let $\phi >0$, $\phi \in C^2(\Omega)$ and $u=\phi v$, then we
have
\begin{equation} \label{eq2}
\int_{\Omega}|\nabla u|^2dx = -\int_{\Omega}\frac{\Delta
\phi}{\phi}u^2dx + \int_{\Omega} \phi^2 |\nabla v|^2dx, ~~~~
\forall u \in C^{\infty}_0(\Omega).
\end{equation}
\end{lemma}
\begin{proof}
By expanding the square we have
\[
\int_{\Omega} |\nabla u+\mathbf{F}u|^2dx =\int_{\Omega}|\nabla
u|^2 dx +\int_{\Omega}|\mathbf{F}|^2u^2dx+\int_{\Omega}\mathbf{F}
\cdot \nabla u^2dx.
\]
Identity (\ref{eq1}) now follows by integrating by parts the
last term.

To prove (\ref{eq2}) we apply (\ref{eq1}) to
$\mathbf{F}=-\frac{\nabla \phi}{\phi}$. Elementary calculations
now yield the result.

\end{proof}

We especially want to study inequalities of the type
\[
  \int_{\mathbb{R}^n_+}|\nabla u|^2dx \geq \int_{\mathbb{R}^n_+}
\left(\frac{\beta_1}{x_1^2}+
 \frac{\beta_2}{x_1^2+x_2^2} + \ldots + \frac{\beta_n}{x_1^2+x_2^2+
\ldots+x_n^2} \right)|u|^2dx,
  \quad  u \in C_0^{\infty}(\mathbb{R}^n_+),
\]
where
 $\mathbf{\beta} = (\beta_1,\ldots,\beta_n)$ is a vector of
 nonnegative constants.
The case when $\beta_1=\frac{1}{4}$ is especially interesting
since it corresponds to the term in the standard Hardy inequality.
So every legitimate choice of $\beta$ with $\beta_1=\frac{1}{4}$
corresponds to an improved Hardy
inequality.
\quad
Let us introduce some notation. Let
\[
\mathbf{X_k} := (x_1,\ldots,x_k,0,\ldots,0) \quad \textrm{so~that }
\quad   |  \mathbf{X_k} |^2  =x_1^2 + \ldots + x_k^2.
\]
We now give the proof of the first part of Theorem A:\\

\noindent
{\em Proof of Theorem A part (i):}
Let  $\gamma_1$, $\gamma_2$, $\ldots$, $\gamma_n$ be arbitrary real numbers
and set
\[
 \phi:=|\mathbf{X_1}|^{-\gamma_1}|\mathbf{X_2}|^{-\gamma_2}\cdot \ldots \cdot
 |\mathbf{X_n}|^{-\gamma_n},
\]
and
\[
 \mathbf{F}: = - \frac{\nabla \phi}{\phi}.
\]
An easy  calculation shows that
\[
 \mathbf{F} = \sum_{m=1}^n  \gamma_m   \frac{\mathbf{X_m}  }
{ |\mathbf{X_m}|^2}.
\]
With this choice of $\mathbf{F}$, we get
\[
\mathrm{div} \mathbf{F} = \sum_{m=1}^n   \gamma_m
  \frac{(m-2)}{ |\mathbf{X_m}|^2 },
\]
and
\[
|\mathbf{F} |^2 = \sum_{m=1}^n \frac{\gamma_m^2}{ |\mathbf{X_m}|^2 }
  + 2 \sum_{m=1}^n  \sum_{j=1}^{m-1}
\gamma_m \gamma_j  \frac{\mathbf{X_m}  }
{ |\mathbf{X_m}|^2}  \frac{\mathbf{X_j}  }
{ |\mathbf{X_j}|^2}
=
 \sum_{m=1}^n
\frac{\gamma_m^2}{ |\mathbf{X_m}|^2 }
+
 2 \sum_{m=1}^n  \sum_{j=1}^{m-1}
 \frac{\gamma_m \gamma_j}
 {|\mathbf{X_j}|^2  }.
\]
We then get that
\be\la{2.5}
-\frac{\Delta \phi}{\phi} =
\mathrm{div} \mathbf{F} - |\mathbf{F}|^2 = \sum_{m=1}^n
\frac{\beta_m}{|\mathbf{X_m}|^2},
\ee
where
\begin{eqnarray*}
 \beta_1 &=& -\gamma_1(\gamma_1+1),  \\
 \beta_m  &=&  - \gamma_m( 2-m +  \gamma_m  + 2 \sum_{j=1}^{m-1} \gamma_j ),
~~~~m=2,3,\ldots,n.
\end{eqnarray*}
We  next  set
\begin{eqnarray*}
 \gamma_1 &=& \alpha_1- \frac12,  \\
 \gamma_m  &=& \alpha_m- \alpha_{m-1} + \frac12,
~~~~m=2,3,\ldots,n.
\end{eqnarray*}
With this choice of $\gamma$'s
 the $\beta$'s are given as in the statement of the Theorem.

As a consequence of Lemma \ref{mainthm} we have that
\be\la{2.6}
 \int_{\mathbb{R}^n_+}|\nabla u|^2 dx \geq
 \int_{\mathbb{R}^n_+}\left(\mathrm{div}\mathbf{F}
  - |\mathbf{F} |^2\right) u^2  dx.
\ee
The result then  follows from (\ref{2.5}) and (\ref{2.6}).

\finedim

\noindent
{\bf Remark} It is easy to check that for any choice of $n$ real numbers
$\xa_1, \ldots, \xa_n$, we can find $n$ nonpositive real numbers
$\xa'_1, \ldots, \xa'_n$ such that they give the same constants
$\xb_1,\ldots, \xb_n$. Consequently, without loss of generality,
we may assume that the real numbers $\xa_1, \ldots, \xa_n$ are
nonpositive.

In the  above  theorem  we have   a lot of freedom. We can choose
the $\gamma$'s in many different ways, each choice giving  a
different inequality. We may, for instance, first maximize
$\beta_1$ and then $\beta_2$ and so on. More generally, we might
try to make the first $m-1$ $\beta_m's$ equal to zero and then
maximize the $\beta_m$'s in increasing order.

In fact we have the following  corollary
\begin{corollary} \label{mainineq}
Let k=1,\ldots,n, then
\begin{eqnarray*}
   \int_{\mathbb{R}^n_+}|\nabla u|^2dx &\geq& \int_{\mathbb{R}^n_+}\left(
  \frac{k^2}{4}\frac{1}{x_1^2+\ldots+x_k^2} + \frac{1}{4}\frac{1}{x_1^2+\ldots+x_{k+1}^2} + \ldots  \right. \\
  &+&  \left. \frac{1}{4}\frac{1}{x_1^2+x_2^2+\ldots+x_n^2} \right)u^2dx,
  \quad  u \in C_0^{\infty}(\mathbb{R}^n_+).
\end{eqnarray*}
\end{corollary}
\begin{proof}
In the case $k=1$ we  choose $\alpha_1= \alpha_2=\ldots = \alpha_n =0$.
 In this case all $\beta_k$'s are equal to $1/4$.

In the general case $k>1$ we  choose $\alpha_m=-m/2$, when  $m=1,2,\ldots,k-1$
and  $\alpha_m=0$,  when  $m=k,\ldots,n$.

\end{proof}

We next  give the proof of the second part of Theorem A: \\

\noindent
{\em Proof of Theorem A, part (ii):}
We will  first  prove that $\beta_1 \leq \frac14$, therefore
  $\beta_1= -\alpha_1^2 +\frac14$,
for suitable  $\alpha_1 \leq 0$. Then, for this $\beta_1$, we will  prove that
$\beta_2 \leq (\alpha_1-\frac12)^2$,   and therefore $\beta_2 = -\alpha_2^2 +
(\alpha_1-\frac12)^2$ for suitable $\xa_2 \leq 0$  and so on.

\noindent
{\bf Step 1.}
Let us first prove  the estimate for  $\beta_1$.  To this end we set
\be\la{22.1}
Q_1[u]:=\frac{\int_{\mathbb{R}^n_+}  |\nabla u|^2dx - \sum_{i=2}^{n}
 \beta_i \int_{\mathbb{R}^n_+}
\frac{u^2}{(x_1^2+x_2^2+\ldots+x_i^2)} dx }{\int_{\mathbb{R}^n_+}
\frac{u^2}{x_1^2} dx }.
\ee
We clearly  have that $\xb_1 \leq \inf_{u \in  C_0^{\infty}(\mathbb{R}^n_+) } Q_1[u]$.
In the sequel we will show that
\be\la{22.2}
\inf_{u \in  C_0^{\infty}(\mathbb{R}^n_+) } Q_1[u]
  \leq  \frac14,
\ee
whence,  $\beta_1 \leq \frac14$.

At this point we introduce a family of cutoff functions for later use.
For $j = 1, \ldots, n$ and $k_j > 0$ we set
\[
\phi_j(t) = \left\{ \begin{array}{ll}
0, &  ~~ t < \frac{1}{k_j^2}\\
1+\frac{\ln{k_jt}}{\ln{k_j}}, & ~~ \frac{1}{k_j^2} \leq t < \frac{1}{k_j}\\
1, & ~~ t \geq \frac{1}{k_j},\\
\end{array} \right.
\]
and
\[
h_{k_j}(x) := \phi_j(r_j) \quad \textrm{ where } \quad
r_j:=|\mathbf{X_j}|=  (x_1^2+\ldots+x_j^2)^{\frac{1}{2}}.
\]
Note that
\[
|\nabla h_{k_j}(x)|^2 = \left\{ \begin{array}{ll}
\frac{1}{\ln^2{k_j}}\frac{1}{r_j^2} &  \frac{1}{k_j^2}
 \leq r_j \leq \frac{1}{k_j}\\
0 &  \textrm{otherwise}\\
\end{array} \right..
\]
We also denote by  $ \phi(x)$ a  radially symmetric
 $C_0^{\infty}(\mathbb{R}^n)$ function
such that  $\phi =1$ for $|x|<1/2$
and $\phi =0$ for $|x|>1$.

To prove (\ref{22.2}) we consider the family of functions
\be\la{22.3}
u_{k_1}(x) = x_1^{\frac12} h_{k_1}(x) \phi(x).
\ee
We will show that as $k_1 \ra \infty$
\be\la{22.4}
\frac{\int_{\mathbb{R}^n_+}  |\nabla u_{k_1}|^2dx -
 \sum_{i=2}^{n} \beta_i \int_{\mathbb{R}^n_+}
\frac{u^2_{k_1}}{(x_1^2+x_2^2+\ldots+x_i^2)} dx }{\int_{\mathbb{R}^n_+}
\frac{u^2_{k_1}}{x_1^2} dx } =
\frac{\int_{\mathbb{R}^n_+}  |\nabla u_{k_1}|^2dx}{\int_{\mathbb{R}^n_+}
\frac{u^2_{k_1}}{x_1^2} dx }   + o(1).
\ee
To see this, let us first examine the behavior of the  denominator.
For $k_1$ large we easily compute
\be\la{22.4a}
\int_{\mathbb{R}^n_+}  \frac{u^2_{k_1}}{x_1^2} dx  =
\int_{\mathbb{R}^n_+}   x_1^{-1}   h_{k_1}^2 \phi^2    dx > C
 \int_{\frac{1}{k_1}}^{\frac12}  x_1^{-1}   dx_1 > C~\ln{k_1}.
\ee
On the other hand by Lebesgue dominated theorem the terms
$ \sum_{i=2}^{n}
 \beta_i \int_{\mathbb{R}^n_+}
\frac{u^2_{k_1}}{(x_1^2+x_2^2+\ldots+x_i^2)} dx $ are easily seen
to be bounded as $k_1 \ra \infty$.  From this and (\ref{22.4a})
we conclude (\ref{22.4}).

We now estimate the gradient term in (\ref{22.4}).
\be\la{22.4b}
\int_{\mathbb{R}^n_+}  |\nabla u_{k_1}|^2 dx =
  \frac14  \int_{\mathbb{R}^n_+}   x_1^{-1}   h_{k_1}^2 \phi^2    dx
+ \int_{\mathbb{R}^n_+}   x_1  | \nabla  h_{k_1}|^{2} \phi^2
+ \int_{\mathbb{R}^n_+}   x_1   h_{k_1}^{2} | \nabla \phi|^{2} +  mixed~~ terms.
\ee
The first integral of the right hand side
 behaves exactly as the denominator, cf (\ref{22.4a}), that is,
 it goes to infinity
like  $O(\ln{k_1})$. The last
integral  is easily seen to be bounded as  $k_1 \ra \infty$. For the
middle integral  we have
\[
 \int_{\mathbb{R}^n_+}   x_1  | \nabla  h_{k_1}|^{2} \phi^2
\leq \frac{C}{\ln^2{k_1}}
 \int_{\frac{1}{k_1^2} \leq x_1 \leq \frac{1}{k_1}}   x_1^{-1} dx_1
\leq \frac{C}{\ln{k_1}}.
\]
As a consequence of these estimates, we easily get that the mixed terms
 in (\ref{22.4b})
are of the order $o(\ln{k_1})$ as  $k_1 \ra \infty$. Hence, we have that
as  $k_1 \ra \infty$,
\be\la{22.4c}
\int_{\mathbb{R}^n_+}  |\nabla u_{k_1}|^2 dx =
 \frac14  \int_{\mathbb{R}^n_+}   x_1^{-1}   h_{k_1}^2 \phi^2    dx
 + o(\ln{k_1}).
\ee
>From (\ref{22.4}), (\ref{22.4a}) and (\ref{22.4c})
 we conclude that  as  $k_1 \ra \infty$
\[
  Q_1[u_{k_1}] = \frac14  + o(1),
\]
hence
$\inf_{u \in  C_0^{\infty}(\mathbb{R}^n_+) } Q_1[u]
  \leq  \frac14$  and consequently  $\xb_1 \leq \frac14$.
Therefore for a suitable nonnegative constant $\alpha_1$ we have that
$\beta_1  =  -\alpha_1^2 + \frac14$. We also set
\be\la{22.5}
\gamma_1 := \alpha_1 -\frac12.
\ee

\noindent
{\bf Step 2.}
We will next show that $\beta_2 \leq  (\alpha_1-\frac12)^2$. To this end,
setting
\be\la{22.6}
Q_2[u] := \frac{\int_{\mathbb{R}^n_+}  |\nabla u|^2dx - (\frac14-\alpha_1^2)
\int_{\mathbb{R}^n_+}
\frac{u^2}{x_1^2} dx -
\sum_{i=3}^{n} \beta_i \int_{\mathbb{R}^n_+}
\frac{u^2}{ |\mathbf{X_i}|^2  } dx }{\int_{\mathbb{R}^n_+}
\frac{u^2 }{ |\mathbf{X_2}|^2 }dx  },
\ee
will  prove that
\[
\inf_{u \in  C_0^{\infty}(\mathbb{R}^n_+) } Q_2[u]
  \leq    (\alpha_1-\frac12)^2.
\]
We now consider the family of functions
\bea\la{22.7}
u_{k_1, k_2}(x) & :=  &
x_1^{-\xg_1} |\mathbf{X_2}|^{\alpha_1-\frac12}
 h_{k_1}(x)h_{k_2}(x) \phi(x)  \nonumber \\
&  =: & x_1^{-\xg_1 }v_{k_1, k_2}(x).
\eea
An a easy calculation shows that
\be\la{22.8}
Q_2[u_{k_1,k_2}]=
\frac{\int_{\mathbb{R}^n_+}x_1^{-2 \xg_1}  |\nabla v_{k_1, k_2} |^2dx -
\sum_{i=3}^{n} \beta_i \int_{\mathbb{R}^n_+}
 x_1^{-2 \xg_1}|\mathbf{X_i}|^{-2}  v_{k_1, k_2}^2  dx }
{\int_{\mathbb{R}^n_+}
x_1^{- 2 \xg_1} |\mathbf{X_2}|^{-2}v_{k_1, k_2}^2  dx  }.
\ee
We  next use the precise form of $v_{k_1,k_2}(x)$. Concerning the
denominator
of $Q_2[u_{k_1,k_2}]$ we have that
\[
\int_{\mathbb{R}^n_+}
x_1^{- 2 \xg_1} |\mathbf{X_2}|^{-2}v_{k_1, k_2}^2  dx
=
\int_{\mathbb{R}^n_+} x_1^{1-2\alpha_1}
(x_1^2+x_2^2)^{\alpha_1-\frac32} h_{k_1}^2 h_{k_2}^2  \phi^2 dx,
\]
Sending $k_1$ to infinity, using the structure of the cutoff functions and
then introducing polar coordinates we  get
\bea\la{22.8a}
\int_{\mathbb{R}^n_+}
x_1^{- 2 \xg_1} |\mathbf{X_2}|^{-2}v_{\infty, k_2}^2  dx &  = &
\int_{\mathbb{R}^n_+} x_1^{1-2\alpha_1}
(x_1^2+x_2^2)^{\alpha_1-\frac32}  h_{k_2}^2  \phi^2 dx \nonumber \\
& \geq &
C \int_{\frac{1}{k_2} < x_1^2 + x_2^2 < \frac12}  x_1^{1-2\alpha_1}
(x_1^2+x_2^2)^{\alpha_1-\frac32} dx_1 dx_2  \\
& \geq &  C     \int_{0}^{\pi}    \int_{\frac{1}{k_2}}^{\frac12} r^{-1}\
 (\sin \theta)^{1- 2\xa_1} dr d \theta \nonumber \\
& \geq &  C  \ln{k_2}. \nonumber
\eea

The terms  in the numerator   that  are multiplied by the  $\beta_i$'s stay
bounded as $k_1$ or $k_2$ go to infinity;  cf the estimates related to
(\ref{es3}) in step 3.
\bea\la{22.8b}
\int_{\mathbb{R}^n_+} x_1^{-2 \xg_1}  |\nabla v_{k_1, k_2}|^2 dx &  = &
  \left(\xa_1-\frac12 \right)^2 \int_{\mathbb{R}^n_+}
  x_1^{-2 \xg_1}  |\mathbf{X_2}|^{2 \xa_1 -3}   h_{k_1}^2  h_{k_2}^2 \phi^2  dx
\nonumber \\
&  & ~+
 \int_{\mathbb{R}^n_+}  x_1^{-2 \xg_1}  |\mathbf{X_2}|^{2 \xa_1 -1}
   |\nabla (h_{k_1} h_{k_2})|^2 \phi^2
 \\
&  & ~+ \int_{\mathbb{R}^n_+}  x_1^{-2 \xg_1}
  |\mathbf{X_2}|^{2 \xa_1 -1}
  h_{k_1}^2  h_{k_2}^2  |\nabla \phi|^{2}
\nonumber \\
&  & ~+ mixed ~ terms
    \nonumber
\eea
The first integral in the right hand side above, is the same as the
denominator of $Q_2$, and therefore is finite as $k_1 \ra \infty$ and increases
like $\ln{k_2}$ as $k_2 \ra \infty$, cf (\ref{22.8a}).
 The last integral is bounded, no matter
how big the $k_1$ and $k_2$ are. Concerning the middle term
we have
\bea\la{22.8c}
M[v_{k_1,k_2}] & := &
 \int_{\mathbb{R}^n_+}  x_1^{-2 \xg_1}  |\mathbf{X_2}|^{2 \xa_1 -1}
   |\nabla (h_{k_1} h_{k_2})|^2 \phi^2 dx  \nonumber \\
& = &  \int_{\mathbb{R}^n_+}  x_1^{-2 \xg_1}  |\mathbf{X_2}|^{2 \xa_1 -1}
   |\nabla h_{k_1}|^2  h_{k_2}^2 \phi^2 dx  +
\int_{\mathbb{R}^n_+}  x_1^{-2 \xg_1}  |\mathbf{X_2}|^{2 \xa_1 -1}
    h_{k_1}^2 |\nabla  h_{k_2}|^2 \phi^2 dx  +    ~mixed~~ term  \nonumber \\
& =:&  ~I_1 + I_2 +mixed~~ term.
\eea
Since
\[
 |\mathbf{X_2}|^{2 \xa_1 -1}  h_{k_2}^2 = r_2^{2 \xa_1 -1} \phi_2(r_2)
\leq C_{k_2},~~~~~~~~~~~ 0< r_2 <1,
\]
we easily get
\[
I_1 \leq \frac{C}{(\ln{k_1})^2} \int_{\frac{1}{k_1^2}}^{\frac{1}{k_1}}
x_1^{-1-2 \xa_1} dx_1,
\]
and therefore, since $\xa_1 \leq 0$, \be\la{22.9} I_1  \leq
\frac{C}{\ln{k_1}}, ~~~~~~~~~~~~~k_1 \ra \infty. \ee Also, since
$h_{k_1}^2  \leq 1$,  we  similarly   get (for any $k_1$)
\be\la{22.9a} I_2 \leq \frac{C}{(\ln{k_2})^2}
\int_{\frac{1}{k_2^2}}^{\frac{1}{k_2}} r_2^{-1} dr_2 ~ \leq
\frac{C}{\ln{k_2}}, ~~~~~~~~~~k_2 \ra \infty. \ee
From
(\ref{22.8c})-- (\ref{22.9a}) we have that as $k_2 \ra \infty$,
\[
M[v_{\infty,k_2}] = o(1).
\]
Returning to (\ref{22.8b}) we have that  as $k_2 \ra \infty$,
\be\la{22.10}
\int_{\mathbb{R}^n_+} x_1^{-2 \xg_1}  |\nabla v_{\infty, k_2}|^2 dx
=   \left(\xa_1-\frac12 \right)^2
\int_{\mathbb{R}^n_+}
x_1^{- 2 \xg_1} |\mathbf{X_2}|^{-2}v_{\infty, k_2}^2  dx
+  o( \ln{k_2}).
\ee
We then have that as  $k_2 \ra \infty$,
\be\la{22.11}
Q_2[u_{\infty,k_2}]=  \left(\xa_1-\frac12 \right)^2 + o(1),
\ee
consequently,  $\beta_2 \leq \left(\xa_1-\frac12 \right)^2$, and therefore
 $\beta_2 = -\alpha_2^2 +
(\alpha_1-\frac12)^2$ for suitable $\xa_2 \leq 0$. We also set
\[
\gamma_2  = \alpha_2- \alpha_{1} + \frac12.
\]

\noindent
{\bf Step 3.} The general case.
At the ($q-1$)th step,  $1 \leq q \leq n$, we have already established
that
\begin{eqnarray*}
 \beta_1 &=& -\alpha_1^2   +\frac14,  \\
 \beta_m  &=&  - \alpha_m^2 + \left(\alpha_{m-1}-\frac12 \right)^2,
~~~~m=2,3,\ldots,q-1,
\end{eqnarray*}
for suitable nonpositive constants $a_i$. Also, we have defined
\begin{eqnarray*}
 \gamma_1 &=& \alpha_1- \frac12,  \\
 \gamma_m  &=& \alpha_m- \alpha_{m-1} + \frac12,
~~~~m=2,3,\ldots,q-1.
\end{eqnarray*}
Our goal for the rest of the proof  is to show  that $\beta_q \leq
\left(\alpha_{q-1}-\frac12 \right)^2$. To this end
  we consider  the quotient
\be
Q_q[u] := \frac{\int_{\mathbb{R}^n_+}  |\nabla u|^2dx
-\sum_{q \neq i=1}^{n} \beta_i \int_{\mathbb{R}^n_+}
\frac{u^2}{ |\mathbf{X_i}|^2  } dx}
{\int_{\mathbb{R}^n_+}
\frac{u^2 }{ |\mathbf{X_q}|^2 }dx  }.
\ee
The test function is now given by
\bea
u_{k_1,k_q}(x) & :=  &
x_1^{-\xg_1}   |\mathbf{X_2}|^{-\xg_2}  \ldots
  |\mathbf{X_{q-1}}|^{-\xg_{q-1}}
   |\mathbf{X_q}|^{\alpha_{q-1}-\frac12}
 h_{k_1}(x)
 h_{k_q}(x) \phi(x)  \nonumber \\
&  =: & x_1^{-\xg_1 }  |\mathbf{X_2}|^{-\xg_2}
 \ldots   |\mathbf{X_{q-1}}|^{-\xg_{q-1}}
    v_{k_q}(x).
\eea
A straightforward calculation shows that
\be\la{22.8v}
Q_q[u_{k_1,k_q}]=
\frac{\int_{\mathbb{R}^n_+}
\prod_{j=1}^{q-1}  |\mathbf{X_{j}}|^{-2 \xg_{j}}
 |\nabla v_{k_1, k_q} |^2dx -
\sum_{i=q+1}^{n} \beta_i \int_{\mathbb{R}^n_+}
\prod_{j=1}^{q-1}  |\mathbf{X_{j}}|^{-2 \xg_{j}}
|\mathbf{X_i}|^{-2}  v_{k_1, k_q}^2  dx }
{\int_{\mathbb{R}^n_+}
\prod_{j=1}^{q-1}  |\mathbf{X_{j}}|^{-2 \xg_{j}}
 |\mathbf{X_q}|^{-2}v_{k_1, k_q}^2  dx }.
\ee
Let us first  see the denominator,
\[
D_q[u_{k_1,k_q}] :=
\int_{\mathbb{R}^n_+}
\prod_{j=1}^{q-1}  |\mathbf{X_{j}}|^{-2 \xg_{j}}
 |\mathbf{X_q}|^{ 2 \xa_{q-1}-3} h_{k_1}(x)  h_{k_q}(x) \phi(x)   dx.
\]
Sending $k_1 \ra \infty$,   we have that $ h_{k_1} \ra  1$ and therefore
\[
D_q[u_{\infty,k_q}] =
\int_{\mathbb{R}^n_+}
\prod_{j=1}^{q-1}  |\mathbf{X_{j}}|^{-2 \xg_{j}}
 |\mathbf{X_q}|^{ 2 \xa_{q-1}-3} h_{k_q}(x) \phi(x)   dx.
\]
To see that this is finite  we note that with
 $B_R^+ :=\{ x \in \mathbb{R}^n:  |x|<R, ~x_1 \geq 0 \}$
\bea
D_q[u_{\infty ,k_q}] & \leq &
 \int_{B^+_1 \cap \{ \frac{1}{k_q^2} \leq r_q \leq
\frac{1}{k_q} \}}\prod_{j=1}^{q-1}  |\mathbf{X_{j}}|^{-2 \xg_{j}}
 |\mathbf{X_q}|^{ 2 \xa_{q-1}-3}dx \nonumber \\
& \leq & C  \int_{ \{ \frac{1}{k_q^2} \leq r_q \leq
\frac{1}{k_q} \}}\prod_{j=1}^{q-1}  |\mathbf{X_{j}}|^{-2 \xg_{j}}
 |\mathbf{X_q}|^{ 2 \xa_{q-1}-3}dx_1 \ldots dx_q.
\eea
To estimate  this, we introduce polar coordinates
 $(x_1, \ldots, x_q) \ra (r_q, \theta_1, \ldots, \theta_{q-1})$.
\begin{eqnarray*}
x_1 &=& r_q\sin{\theta_{q-1}}\sin{\theta_{q-2}}\cdot \ldots \cdot
\sin{\theta_2} \sin{\theta_1} \\
x_2 &=& r_q\sin{\theta_{q-1}}\sin{\theta_{q-2}}\cdot \ldots \cdot
\sin{\theta_2} \cos{\theta_1}\\
x_3 &=& r_q\sin{\theta_{q-1}}\sin{\theta_{q-2}}\cdot \ldots \cdot
 \cos{\theta_2} \\
 &\vdots& \\
 x_q &=& r_q \cos{\theta_{q-1}},
\end{eqnarray*}
where $0\leq\theta_1<2\pi$ and $0\leq\theta_m<\pi$ for
$m=2, \ldots, q-1$.
The surface measure on the unit sphere  $S^{q-1}$ then becomes
\[
C (\sin{\theta_{q-1}})^{q-2}(\sin{\theta_{q-2}})^{q-3}\cdots
 \sin{\theta_2}d\theta_1\ldots d\theta_{q-1}.
\]
Also,  $r_q= |\mathbf{X_{q}}|$
and for $1 \leq m \leq q-1$,
\[
r_m=  |\mathbf{X_{m}}|=
 (x_1^2+\ldots+x_m^2)^{\frac12} = r_q \,
 \sin{\theta_{q-1}}\sin{\theta_{q-2}}\cdot
 \ldots \cdot \sin{\theta_{m}}.
\]
We then have
\bea
 \int_{ \{ \frac{1}{k_q^2} \leq r_q   \leq
\frac{1}{k_q} \}}
\prod_{j=1}^{q-1}  |\mathbf{X_{j}}|^{-2 \xg_{j}}
 |\mathbf{X_q}|^{ 2 \xa_{q-1}-3}  dx_1 \ldots dx_q
 &  =  & C \int_{\{ \frac{1}{k_q^2} \leq r_q    \leq
\frac{1}{k_q} \}}
 r_q^{-1}  \prod_{j=1}^{q-1} ( \sin{\theta_j} )^{1-2 \xa_j}
 d \theta_1 \ldots d \theta_{q-1}  d r_q  \nonumber \\
& \leq & C \ln{k_q}.
\eea
On the other hand since,
\[
D_q[u_{\infty ,k_q}]  \geq
 \int_{B^+_{1/2} \cap \{ \frac{1}{k_q} \leq r_q \leq
\frac{1}{2} \}}\prod_{j=1}^{q-1}  |\mathbf{X_{j}}|^{-2 \xg_{j}}
 |\mathbf{X_q}|^{ 2 \xa_{q-1}-3}dx,
\]
by practically the same argument we have that as $k_q \ra \infty$,
\be
D_q[u_{\infty ,k_q}]  \geq  C \ln{k_q}.
\ee

For  $i=q+1, \ldots, n$,
we  consider the terms
\bea\la{es3}
 \int_{\mathbb{R}^n_+}
\prod_{j=1}^{q-1}  |\mathbf{X_{j}}|^{-2 \xg_{j}}
|\mathbf{X_i}|^{-2}  v_{k_1, k_q}^2  dx  &  = &
 \int_{\mathbb{R}^n_+}
\prod_{j=1}^{q-1}  |\mathbf{X_{j}}|^{-2 \xg_{j}}
|\mathbf{X_q}|^{2 \xa_{q-1} -1}
|\mathbf{X_i}|^{-2}  h_{k_1}^2  h_{k_q}^2 \phi^2(x) dx
\nonumber \\
& \leq &
\int_{\mathbb{R}^n_+}
\prod_{j=1}^{q-1}  |\mathbf{X_{j}}|^{-2 \xg_{j}}
|\mathbf{X_q}|^{2 \xa_{q-1} -1}
|\mathbf{X_{q+1}}|^{-2}  h_{k_1}^2  h_{k_q}^2 \phi^2(x) dx.
\eea
Taking first the limit  $k_1 \ra \infty$ and then $k_q \ra \infty$,
the above integral converges to
\[
I_q:=
\int_{\mathbb{R}^n_+}
\prod_{j=1}^{q-1}  |\mathbf{X_{j}}|^{-2 \xg_{j}}
|\mathbf{X_q}|^{2 \xa_{q-1} -1}
|\mathbf{X_{q+1}}|^{-2}  \phi^2(x) dx.
\]
To see that this is finite we introduce
 polar coordinates in $(x_1, \ldots, x_{q+1}) \ra (r_{q+1}, \theta_1
 \ldots,\theta_q)$
and  use elementary estimates  to get
\[
I_q \leq C \int_{B_1^+} \sin \theta_q
\prod_{j=1}^q ( \sin{\theta_j} )^{1-2 \xa_j}
 d \theta_1 \ldots d \theta_{q}  d r_{q+1} < \infty.
\]

We next consider the gradient term
\bea\la{22.50}
\int_{\mathbb{R}^n_+}
\prod_{j=1}^{q-1}  |\mathbf{X_{j}}|^{-2 \xg_{j}}
  |\nabla v_{k_1, k_q}|^2 dx &  = &
  \left(\xa_{q-1}-\frac12 \right)^2 \int_{\mathbb{R}^n_+}
  \prod_{j=1}^{q-1}  |\mathbf{X_{j}}|^{-2 \xg_{j}}
       |\mathbf{X_q}|^{2 \xa_{q-1} -3}   h_{k_1}^2  h_{k_q}^2 \phi^2  dx
\nonumber \\
&  & ~+
 \int_{\mathbb{R}^n_+}
\prod_{j=1}^{q-1}  |\mathbf{X_{j}}|^{-2 \xg_{j}}
  |\mathbf{X_q}|^{2 \xa_{q-1} -1}
   |\nabla (h_{k_1} h_{k_q})|^2 \phi^2
 \\
&  & ~+ \int_{\mathbb{R}^n_+}
 \prod_{j=1}^{q-1}  |\mathbf{X_{j}}|^{-2 \xg_{j}}
  |\mathbf{X_q}|^{2 \xa_{q-1} -1}
  h_{k_1}^2  h_{k_q}^2  |\nabla \phi|^{2}
\nonumber \\
&  & ~+ mixed ~ terms
    \nonumber
\eea The first term of the right hand side is the same as the
denominator. Using polar coordinates and arguments similar to  the
ones used in estimating the gradient term in  (\ref{22.8b}),
 all   other terms of (\ref{22.50})  are bounded as $k_1 \ra \infty$ and
$k_q \ra \infty$. In particular we end up with
\[
\int_{\mathbb{R}^n_+}
\prod_{j=1}^{q-1}  |\mathbf{X_{j}}|^{-2 \xg_{j}}
  |\nabla v_{\infty, k_q}|^2 dx =
 \left(\xa_{q-1}-\frac12 \right)^2 D_q[u_{\infty ,k_q}]
+ o (\ln{k_q}), ~~~~~~k_q \ra \infty.
\]
Putting things together  we have that
\[
Q_q[ u_{\infty, k_q}] =  \left(\xa_{q-1}-\frac12 \right)^2 + o(1),
~~~~~~k_q \ra \infty,
\]
from which it follows that $\beta_q \leq \left(\xa_{q-1}-\frac12 \right)^2$.
This completes the proof of the Theorem.

\finedim

The previous analysis can also lead to the following result:

\begin{theorem}
Let  $\xa_1, \ldots, \xa_k$,  $1 \leq k \leq n-1$,
be  nonpositive constants  and
\begin{eqnarray*}
 \beta_1 &=& -\alpha_1^2   +\frac14,  \\
 \beta_m  &=&  - \alpha_m^2 + \left(\alpha_{m-1}-\frac12 \right)^2,
~~~~m=2,3,\ldots,k.
\end{eqnarray*}
Suppose that there exists a constant $\xb_{k+1}$ such that the
 following
inequality holds
\be\la{2.80}
  \int_{\mathbb{R}^n_+}|\nabla u|^2dx \geq
\int_{\mathbb{R}^n_+}\left(\frac{\beta_1}{x_1^2}+
 \frac{\beta_2}{x_1^2+x_2^2} + \ldots
 + \frac{\beta_{k+1}}{x_1^2+x_2^2+\ldots+x_{k+1}^2} \right)u^2dx,
\ee
for any $u \in C_0^{\infty}(\mathbb{R}^n_+)$.  Then
\be\la{2.81}
\xb_{k+1} \leq  \left(\alpha_{k}-\frac12 \right)^2.
\ee
Moreover,
\be\la{2.82}
\inf_{u \in C^{\infty}_0(\mathbb{R}^n_+)} \frac{
\int_{\mathbb{R}^n_+}|\nabla u|^2dx - \xb_1
\int_{\mathbb{R}^n_+}\frac{|u|^2}{x_1^2}dx - \ldots - \xb_k
\int_{\mathbb{R}^n_+}\frac{|u|^2}{x_1^2+\ldots+x_k^2}dx}{\int_{\mathbb{R}^n_+}
\frac{|u|^2}{x_1^2+\ldots+x_{k+1}^2}dx} = \left(\alpha_{k}-\frac12 \right)^2.
\ee
\end{theorem}
\begin{proof}
The proof of the first part, that is, estimate (\ref{2.81}),
 is contained in the proof of Theorem A(ii).

To establish  the second result (\ref{2.82}), we  first use (\ref{2.81})
to obtain that the infimum  in (\ref{2.82})  is less that or equal to
$\left(\alpha_{k}-\frac12 \right)^2$.  To obtain the reverse inequality
we use Theorem A(i) with $a_{k+l}=-\frac{l-1}{2}$,  $l=1, \ldots, n-k$.
For this choice we have that $\xb_{k+2} = \ldots = \xb_n=0$.

\end{proof}

The following is an interesting consequence of the previous Theorem.

\begin{corollary}
For $1 \leq k \leq n$,
\begin{equation} \label{sats1}
\inf_{u \in C^{\infty}_0(\mathbb{R}^n_+)}
\frac{\int_{\mathbb{R}^n_+}|\nabla u|^2dx}{\int_{\mathbb{R}^n_+}
\frac{|u|^2}{x_1^2+\ldots+x_{k}^2}dx} = \frac{k^2}{4},
\end{equation}
 and
\begin{equation} \label{sats2}
\inf_{u \in C^{\infty}_0(\mathbb{R}^n_+)} \frac{
\int_{\mathbb{R}^n_+}|\nabla u|^2dx - \frac{k^2}{4}
\int_{\mathbb{R}^n_+}\frac{|u|^2}{x_1^2+\ldots+x_k^2}dx-
\frac{1}{4}\int_{\mathbb{R}^n_+}\frac{|u|^2}{x_1^2+\ldots+x_{k+1}^2}dx
-\ldots- \frac{1}{4}\int_{\mathbb{R}^n_+}
\frac{|u|^2}{x_1^2+\ldots+x_{m}^2}dx}{\int_{\mathbb{R}^n_+}
\frac{|u|^2}{x_1^2+\ldots+x_{m+1}^2}dx} = \frac{1}{4}
\end{equation}
for $k\leq m<n$.
\end{corollary}
\begin{proof}
To establish (\ref{sats1}) we use  (\ref{2.82})  with  $\xa_l=-\frac{l}{2}$,
$l=1, \ldots,k-1$.

To establish (\ref{sats2}) we again use  (\ref{2.82})  with
$\xa_l=-\frac{l}{2}$,
$l=1, \ldots,k-1$,  and   $\xa_l =0$, $k \leq l \leq m$. With choice
we have that $\xb_1=\ldots \xb_{k-1}=0$, $\xb_k=\frac{k^2}{4}$ and
$\xb_{l}=\frac14$, $l=k-1, \ldots,m$.

\end{proof}

\section{Hardy-Sobolev-Maz'ya inequalities}

We begin by proving    Theorem C.

\noindent
{\em Proof of Theorem C:}
Our starting point is the Gagliardo-Nirenberg-Sobolev inequality
\begin{equation}
 C_n \int_{\Omega}|f|^{\frac{n}{n-1}}dx \leq \left( \int_{\Omega}|\nabla f|dx \right)^{\frac{n}{n-1}},
 \quad \quad f \in C_0^{\infty}(\Omega). \label{nirenberg}
\end{equation}
Let $f=\phi^{\alpha}w$, where $\alpha = \frac{2(n-1)}{n-2}$. This
leads to
\[
 C_n \int_{\Omega}\phi^{\frac{\alpha n}{n-1}}|w|^{\frac{n}{n-1}}dx \leq
 \left( \int_{\Omega}\alpha \phi^{\alpha-1}|\nabla \phi||w|
 + \phi^{\alpha}|\nabla w|dx \right)^{\frac{n}{n-1}},
 \quad \quad w \in C_0^{\infty}(\Omega).
\]
We  now estimate the first term in the integral according to
inequality (\ref{c1}) and let
$w=|v|^{\theta}$. Then we get
\begin{eqnarray*}
 C \left( \int_{\Omega}\phi^{\frac{\alpha n}{n-1}}|v|^{\frac{\theta
 n}{n-1}}dx\right)^{\frac{n-1}{n}}
 &\leq&
 \int_{\Omega} \phi^{\alpha}|v|^{\theta-1}|\nabla v|dx \\ &\leq&
 \left( \int_{\Omega} \phi^{2\alpha-2}|v|^{2\theta-2}
  dx \right)^{1/2} \left(\int_{\Omega} \phi^{2}|\nabla
 v|^2dx \right)^{1/2}
\end{eqnarray*}
 The choice
\[
  \theta = \alpha = \frac{2(n-1)}{n-2}
\]
gives us the inequality
\begin{equation}
  C \left( \int_{\Omega}\phi^{\frac{2n}{n-2}}|v|^{\frac{2n}{n-2}}
  dx\right)^{\frac{n-2}{n}} \leq \int_{\Omega} \phi^{2}|\nabla
 v|^2dx. \label{sobolev}
\end{equation}
Let $u=\phi v$. By lemma (\ref{mainthm}) we have
\[
\int_{\Omega}|\nabla u|^2dx = -\int_{\Omega}\frac{\Delta
\phi}{\phi}u^2dx + \int_{\Omega} \phi^2 |\nabla v|^2dx.
\]
We conclude the proof by combining this result with
(\ref{sobolev}).
\finedim

Condition (\ref{c1}) might seem to be unnatural and not
easily checked. However, it will be very natural and is easily
verified for our choices of $\phi$.

 To produce Hardy inequalities in the half-space with remainder terms also
including the Sobolev term, we will need a weighted version of the
Sobolev inequality.
\begin{theorem} \label{weight}
Let $\xs_1, \xs_2,  \ldots, \xs_k$  be real numbers  for some $k$ with $1\leq k \leq n$.  We set $c_l :=|\xs_1 +\ldots +\xs_l +l-1|$, for  $1 \leq l \leq k$.
We assume that
\[
c_l \neq 0 ~~~~~~~~~{\rm   whenever} ~~~~~~~~~ \xs_l \neq 0.
\]
 Then, there exists a positive constant $C$ such that
 for any $w   \in C_0^{\infty}(\mathbb{R}^n_+)$ there holds
\be\la{3.5}
\int_{\mathbb{R}^n_+}x_1^{\xs_1}
 |\mathbf{X_2}|^{\xs_2}
 \ldots|\mathbf{X_k}|^{\xs_k} |\nabla w|dx
\geq
C \left(
 \int_{\mathbb{R}^n_+}   \left( x_1^{\xs_1}
 |\mathbf{X_2}|^{\xs_2}
 \ldots|\mathbf{X_k}|^{\xs_k} |w| \right)^{\frac{n}{n-1}} dx
 \right)^{\frac{n-1}{n}},
\ee
and
\[
\int_{\mathbb{R}^n_+}   x_1^{\frac{\xs_1 (n-2)}{(n-1)}}
 |\mathbf{X_2}|^{\frac{\xs_2 (n-2)}{(n-1)}}
 \ldots|\mathbf{X_k}|^{\frac{\xs_k (n-2)}{(n-1)}}
 |\nabla w|^2 dx
\geq
\]
\be\la{3.6}
\geq C
 \left( \int_{\mathbb{R}^n_+}   \left( x_1^{\frac{\xs_1 (n-2)}{2(n-1)}}
 |\mathbf{X_2}|^{\frac{\xs_2 (n-2)}{2(n-1)}}
 \ldots|\mathbf{X_k}|^{\frac{\xs_k (n-2)}{2(n-1)}}
  |w| \right)^{\frac{2n}{n-2}} dx \right)^{\frac{n-2}{n}}.
\ee

\end{theorem}

\begin{proof}
For  $\Omega = \mathbb{R}^n_+$  we let  $u=x_1^{\xs_1}v$
in the Sobolev inequality (\ref{nirenberg}) to get
\[
 C_n \int_{\mathbb{R}^n_+}x_1^{\frac{n\xs_1}{n-1}}
|v|^{\frac{n}{n-1}}dx \leq \left( \int_{\mathbb{R}^n_+}
 |\xs_1| x_1^{\xs_1-1}|v| + x_1^{\xs_1}|\nabla v|dx \right)^{\frac{n}{n-1}},
 \quad \quad v \in C_0^{\infty}(\mathbb{R}^n_+).
\]
Using the inequality
\begin{equation}
  \left| \int_{\mathbb{R}^n_+} \mathrm{div} \mathbf{F} |v|dx \right| \leq \int_{\mathbb{R}^n_+}| \mathbf{F}||\nabla v|dx, \label{L1}
\end{equation}
with the vector field $(x_1^{\xs_1},0,\ldots,0)$ one obtains
\[
 |\xs_1| \int_{\mathbb{R}^n_+}x_1^{\xs_1-1}|v| dx \leq \int_{\mathbb{R}^n_+} x_1^{\xs_1}|\nabla v|dx
\]
and hence that
\[
 C_n \int_{\mathbb{R}^n_+}x_1^{\frac{n\xs_1}{n-1}}|v|^{\frac{n}{n-1}}dx \leq \left( \int_{\mathbb{R}^n_+}
  x_1^{\xs_1}|\nabla v|dx \right)^{\frac{n}{n-1}},
 \quad \quad v \in C_0^{\infty}(\mathbb{R}^n_+).
\]
Now let $v=  |\mathbf{X_2}|^{\xs_2}w  = (x_1^2+x_2^2)^{\xs_2/2}w$ in the above inequality.
This gives
\begin{eqnarray*}
   C_n \int_{\mathbb{R}^n_+}x_1^{\frac{n\xs}{n-1}}(x_1^2+x_2^2)^{\frac{n\xs_2}{2(n-1)}}|w|^{\frac{n}{n-1}}dx
  &\leq& \left(  \int_{\mathbb{R}^n_+}
  x_1^{\xs_1}(x_1^2+x_2^2)^{\xs_2/2}|\nabla w|dx \right. \\
  &+& \left.\int_{\mathbb{R}^n_+} | \xs_2| x_1^{\xs_1}(x_1^2+x_2^2)^{\xs_2/2-1/2}  |w|
  \right)^{\frac{n}{n-1}}.
\end{eqnarray*}
Letting  $\mathbf{F} =
x_1^{\xs_1}(x_1^2+x_2^2)^{\xs_2/2-1/2}\mathbf{X}_2    $ in (\ref{L1}),
we  get
\be\la{3.11}
 |\xs_1 + \xs_2 + 1|\int_{\mathbb{R}^n_+}x_1^{\xs_1}(x_1^2+x_2^2)^{\xs_2/2-1/2}|w|dx
 \leq \int_{\mathbb{R}^n_+} x_1^{\xs_1}(x_1^2+x_2^2)^{\xs_2/2}|\nabla w|dx.
\ee
Combining the previous two estimates we conclude
\[
 c\int_{\mathbb{R}^n_+}x_1^{\frac{n\xs_1}{n-1}}(x_1^2+x_2^2)^{\frac{n\xs_2}{2(n-1)}}|w|^{\frac{n}{n-1}}dx
  \leq   \left( \int_{\mathbb{R}^n_+}
  x_1^{\xs_1}(x_1^2+x_2^2)^{\xs_2/2}|\nabla w|dx
  \right)^{\frac{n}{n-1}}.
\]
Note that, in case $\xs_2=0$, we  have the desired result
immediately and we  do  not  have to check  whether  the constant $\xs_1 +
\xs_2 + 1$ is zero or not.
\quad
 We may repeat this procedure iteratively. In the $l$-th
step we  need the analogue of (\ref{3.11}) which is
\begin{eqnarray*}
 &c_l& \int_{\mathbb{R}^n_+}x_1^{\xs_1}(x_1^2+x_2^2)^{\frac{\xs_2}{2}} \cdot \ldots \cdot
  (x_1^2+\ldots+x_l^2)^{\frac{\xs_l-1}{2}}|w|dx \\
 &\leq& \int_{\mathbb{R}^n_+} x_1^{\xs_1}(x_1^2+x_2^2)^{\frac{\xs_2}{2}}
   \cdot \ldots \cdot (x_1^2+\ldots+x_l^2)^{\frac{\xs_l}{2}} |\nabla w|dx
\end{eqnarray*}
for some positive constant $c_l$.
This  follows from (\ref{L1})  with
\[
\mathbf{F} =
x_1^{\xs_1}(x_1^2+x_2^2)^{\frac{\xs_2}{2}}    \ldots
  (x_1^2+\ldots+x_l^2)^{\frac{\xs_l-1}{2}} \mathbf{X}_l,
\]
there.
 For  this choice we get
\[
  c_l = |\xs_1 + \ldots + \xs_l + (l-1)|.
\]
So our procedure works nicely in case $c_l \neq 0$ for those $l$
such that $\xs_l \neq 0$. This proves (\ref{3.5}).

To show (\ref{3.6}) we apply (\ref{3.5}) to the function
 $w=|v|^{\theta}$.
Trivial estimates give
\begin{eqnarray*}
 &C&\int_{\Omega}x_1^{\frac{n\xs_1}{n-1}}(x_1^2+x_2^2)^{\frac{n\xs_2}{2(n-1)}} \cdot \ldots \cdot
  (x_1^2+\ldots+x_k^2)^{\frac{n\xs_k}{2(n-1)}}|v|^{\frac{n \theta}{n-1}}dx \\
 &\leq& \left( \theta \int_{\Omega} x_1^{\xs_1}(x_1^2+x_2^2)^{\frac{\xs_2}{2}}
   \cdot \ldots \cdot (x_1^2+\ldots+x_k^2)^{\frac{\xs_k}{2}} |v|^{\theta-1}|\nabla v|dx
   \right)^{\frac{n}{n-1}}.
\end{eqnarray*}
 We will then apply H\"olders inequality to the
 right hand side. We want to do it in such a way that
 one of the factors becomes
 identical to the left hand side raised to some power. Therefore
 we need to choose $\theta$ so that
 \[
 \frac{n \theta}{n-1}=2\theta -2 \quad \Leftrightarrow \quad \theta =
 \frac{2(n-1)}{n-2}.
 \]
H\"olders inequality then immediately gives the result.

\end{proof}

We are now ready to give the proof of Theorem B: \\

\noindent
{\em  Proof of Theorem B:} For $\phi>0$ and $u= \phi v$,
  Lemma  \ref{mainthm} gives us the inequality
\be\la{3.30}
\int_{\mathbb{R}^n_+}|\nabla u|^2 dx + \int_{\mathbb{R}^n_+}\frac{\Delta
\phi}{\phi}|u|^2dx
 \geq  \int_{\mathbb{R}^n_+} \phi^2|\nabla v|^2dx.
\ee
We will choose for $\phi$,
\bea
\phi(x) &  =  &  \left( x_1^{\xs_1}\cdot
  (x_1^2+x_2^2)^{\frac{\xs_2}{2}}\cdot
  \ldots \cdot (x_1^2+\ldots +
  x_n^2)^{\frac{\xs_n}{2}}\right)^{\frac{n-2}{2(n-1)}}  \nonumber \\
 &  =  &
|\mathbf{X_1}|^{-\gamma_1}|\mathbf{X_2}|^{-\gamma_2}\cdot \ldots \cdot
 |\mathbf{X_n}|^{-\gamma_n},
\eea
where,
\begin{eqnarray*}
 \gamma_1 &=& \alpha_1- \frac12,  \\
 \gamma_m  &=& \alpha_m- \alpha_{m-1} + \frac12,
~~~~m=2,3,\ldots,n.
\end{eqnarray*}
and
\[
\xs_m = - \frac{2(n-1)}{n-2}\gamma_{m}~~~~~~~~~~m=1, \ldots,n.
\]
We now apply (\ref{3.6}) of Theorem \ref{weight} to obtain that
\[
 \int_{\mathbb{R}^n_+} \phi^2|\nabla v|^2dx
\geq C
\left( \int_{\mathbb{R}^n_+}
|\phi v|^{\frac{2n}{n-2}}dx \right)^{\frac{n-2}{n}},
\]
provided that
\be\la{3.19}
c_l :=|\xs_1 +\ldots +\xs_l +l-1| \neq 0, ~~~~~
whenever ~~~\xs_l \neq 0,
\ee
 for  $1 \leq l \leq n$.
Combining this with (\ref{3.30}) we get
\[
\int_{\mathbb{R}^n_+}|\nabla u|^2 dx + \int_{\mathbb{R}^n_+}\frac{\Delta
\phi}{\phi}|u|^2dx
 \geq C \left( \int_{\mathbb{R}^n_+}
|u|^{\frac{2n}{n-2}}dx \right)^{\frac{n-2}{n}}
\]
On the other hand,  by Theorem A(i),
\[
-\frac{\Delta
\phi}{\phi}
=
\frac{\beta_1}{x_1^2}+
 \frac{\beta_2}{x_1^2+x_2^2} + \ldots +
 \frac{\beta_n}{x_1^2+x_2^2+\ldots+x_n^2},
\]
and the desired inequality follows. It remains to check condition (\ref{3.19}).
After some elementary calculations we see that
\[
c_l = \frac{2(n-1)}{n-2} \left|\xa_l - \frac{n-l}{2(n-1)} \right|, ~~~~~~~
l=1,\ldots, n.
\]
Since  $\xa_l \leq 0$  we clearly have that $c_l \neq 0$ for $l=1,\ldots, n-1$.
Moreover $c_n \neq 0$ when $\xa_n < 0$.  This completes the proof of (\ref{1.20}).

In the rest of the proof we will show that (\ref{1.20}) fails in case $\xa_n=0$.
To this end we will establish that
\be\la{3.25}
\inf_{u \in C^{\infty}_0(\mathbb{R}^n_+)} \frac{
\int_{\mathbb{R}^n_+}|\nabla u|^2dx - \xb_1
\int_{\mathbb{R}^n_+}\frac{|u|^2}{x_1^2}dx - \ldots - \xb_n
\int_{\mathbb{R}^n_+}\frac{|u|^2}{x_1^2+\ldots+x_n^2}dx}{ \left(
\int_{\mathbb{R}^n_+}|u|^{\frac{2n}{n-2}}dx\right)^{\frac{n-2}{n}}}=0,
\ee
where  $\xb_n =  \left(\alpha_{n-1}-\frac12 \right)^2$. Let
\[
u(x)=
 x_1^{-\xg_1 }  |\mathbf{X_2}|^{-\xg_2}
 \ldots   |\mathbf{X_{n-1}}|^{-\xg_{n-1}}
    v(x).
\]
A straightforward calculation, quite similar to the one leading to (\ref{22.8}), shows
that  the infimum in (\ref{3.25}) is the same as the following  infimum
\be\la{3.60}
\inf_{v \in C^{\infty}_0(\mathbb{R}^n_+)}
\frac{\int_{\mathbb{R}^n_+}
\prod_{j=1}^{n-1}  |\mathbf{X_{j}}|^{-2 \xg_{j}}
 |\nabla v |^2dx -
 \beta_n \int_{\mathbb{R}^n_+}
\prod_{j=1}^{n-1}  |\mathbf{X_{j}}|^{-2 \xg_{j}}
|\mathbf{X_n}|^{-2}  v ^2  dx }
{\left( \int_{\mathbb{R}^n_+}
\left(\prod_{j=1}^{n-1}  |\mathbf{X_{j}}|^{- \xg_{j}} \right)^{\frac{2n}{n-2}}
  |v|^{\frac{2n}{n-2}}  dx \right)^{\frac{n-2}{n}}}.
\ee
We  now choose the following test functions
\be\la{tf}
v_{k_1, \xe}=  |\mathbf{X_{n}}|^{- \xg_{n} +\xe} h_{k_1}(x) \phi(x),~~~~~~~~~~~~~\xe>0,
\ee
where $ h_{k_1}(x)$ and $\phi(x)$ are the same test functions as in  the first step
of the proof of Theorem  A(ii). For this choice, after straightforward
calculations,  quite similar to the ones used in the proof of Theorem A(ii),
we obtain the following estimate for the numerator $N$ in (\ref{3.60}).
\begin{eqnarray*}
N[v_{\infty, \xe}] &  =  & \left( \left(\alpha_{n-1}-\frac12  +\xe \right)^2-
\left(\alpha_{n-1}-\frac12  \right)^2 \right)
\int_{\mathbb{R}^n_+}
\prod_{j=1}^{n-1}  |\mathbf{X_{j}}|^{-2 \xg_{j}}
 |\mathbf{X_{n}}|^{- 2 \xg_{n}+2 +\xe} \phi^2(x) dx + O_{\xe}(1),  \\
& =  &
 C  \xe  \int_{\mathbb{R}^n_+}
 r^{-1+ 2 \xe}  \prod_{j=1}^{n} ( \sin{\theta_j} )^{1-2 \xa_j}  \phi^2(r)
 d \theta_1 \ldots d \theta_{n-1}  d r + O_{\xe}(1)    \\
& =  & C \xe  \int_{0}^{1} r^{-1+\xe} dr + O_{\xe}(1).
\end{eqnarray*}
In the above calculations we have taken the limit $k_1 \ra \infty$
and we have used polar coordinates in $(x_1,\ldots, x_n) \ra
(\theta_1,\ldots, \theta_{n-1},r)$. We then conclude that
\be\la{3.81} N[v_{\infty, \xe}]  < C, ~~~~~~~~~~~{\rm as}
~~~~~~~~\xe \ra 0. \ee Similar calculations for the denominator
$D$  in (\ref{3.60}) reveal that
\begin{eqnarray*}
D[v_{\infty, \xe}] &  =  & C
\left(\int_{\mathbb{R}^n_+}  r^{-1+\frac{2  \xe n}{n-2}}
\prod_{j=1}^{n-1}(\sin{\theta_j})^{\frac{n-j}{n-2}-\frac{2n \xa_j}{n-2} -1 }
\phi^{\frac{2n}{n-2}}
 d\theta_1 \ldots d\theta_{n-1} dr   \right)^{\frac{n-2}{n}}  \\
 &  \geq   & C \left( \int_{0}^{\frac12}
 r^{-1+\frac{2  \xe n}{n-2}} dr \right)^{\frac{n-2}{n}} \\
& = & C \xe^{- \frac{n-2}{n}}.
\end{eqnarray*}
We then  have that
\[
\frac{N[v_{\infty, \xe}]}{D[v_{\infty, \xe}] } \ra 0 ~~~~~~{\rm as}~~~~
\xe \ra 0,
\]
and therefore the  infimum in (\ref{3.60})  or (\ref{3.25}) is equal to zero.
This completes the proof of the Theorem.

\finedim

Here is a consequence of the Theorem  B.

\begin{corollary}
Let  $1 \leq
k<n$.  For any  $\beta_n < \frac{1}{4}$, there exists a positive constant $C$
such that for  all  $u \in C^{\infty}_0(\mathbb{R}^n_+)$ there holds
\begin{eqnarray*}
  \int_{\mathbb{R}^n_+}|\nabla u|^2dx &\geq& \int_{\mathbb{R}^n_+}
 \left( \frac{k^2}{4}\frac{1}{x_1^2+x_2^2+\ldots+x_k^2}
 + \frac{1}{4}\frac{1}{x_1^2+x_2^2+\ldots+x_{k+1}^2}
 +\ldots \right.\\ &+& \left. \frac{1}{4}\frac{1}{x_1^2+x_2^2+\ldots+x_{n-1}^2}+
 \frac{\beta_n}{x_1^2+x_2^2+\ldots+x_{n}^2}
 \right)|u|^2dx
  + C\left(\int_{\mathbb{R}^n_+}|u|^{2^*}dx
  \right)^{\frac{2}{2^*}},
\end{eqnarray*}
If $\beta_n = \frac{1}{4}$ the previous inequality fails.

In case  $k=n$ we have that for  any  $\beta_n < \frac{n^2}{4}$,
 there exists a positive constant $C$
such that for  all  $u \in C^{\infty}_0(\mathbb{R}^n_+)$ there holds
\[
 \int_{\mathbb{R}^n_+}|\nabla u|^2dx \geq
 \xb_n\int_{\mathbb{R}^n_+}\frac{|u|^2}{x_1^2+x_2^2+\ldots+x_{n}^2}dx + C\left(\int_{\mathbb{R}^n_+}|u|^{2^*}dx
  \right)^{\frac{2}{2^*}}.
\]
The above inequality fails for $\beta_n = \frac{n^2}{4}$
\end{corollary}
\begin{proof} In Theorem B we make the following choices:
In the case $k=1$ we  choose $\alpha_1= \alpha_2=\ldots = \alpha_{n-1} =0$.
 In this case  $\beta_k=1/4$,    $k=1,\ldots,n-1$. The condition $\xa_n<0$ is
equivalent to $\beta_n < \frac{1}{4}$.

In the case  $1<k \leq n-1$ we  choose $\alpha_m=-m/2$, when  $m=1,2,\ldots,k-1$
and  $\alpha_m=0$,  when  $m=k,\ldots,n-1$.
Finally, in case  $k=n$,  we  choose $\alpha_m=-m/2$, for  $m=1,2,\ldots,n-1$.

\end{proof}

\section{Further generalizations}

The techniques used in the previous sections can be generalized to
other situations as well. For example, consider the subset of
$\mathbb{R}^n$, where $x_1,x_2, \ldots ,x_k
>0$. We denote  this domain  by $\mathbb{R}^n_{k_+}$.
Then we can easily prove the Hardy-Sobolev inequality
\begin{theorem}
There exists a positive constant $C$ such that for  any
 $u \in C^{\infty}_0(\mathbb{R}^n_{k_+})$  there holds
\[
\int_{\mathbb{R}^n_{k_+}}|\nabla u|^2dx \geq
\frac{1}{4}\int_{\mathbb{R}^n_{k_+}}
 \left(\frac{1}{x_1^2}+\ldots+\frac{1}{x_k^2}\right)|u|^2dx + C\left(
\int_{\mathbb{R}^n_{k_+}}|u|^{2^*}dx \right)^{\frac{2}{2^*}}.
 \]
\end{theorem}
\begin{proof}
Let  $\phi = \sqrt{x_1\cdot \ldots \cdot
x_k}$.  For  $u=\phi w$  we calculate
to get
\begin{eqnarray*}
\int_{\mathbb{R}^n_{k_+}} |\nabla u|^2dx
 &= &\int_{\mathbb{R}^n_{k_+}} |\sqrt{x_1\cdot \ldots \cdot x_k} \cdot \nabla w + \frac{1}{2}
 \sqrt{x_1\cdot \ldots \cdot x_k}
 \left( \frac{1}{x_1}, \ldots,
 \frac{1}{x_k} \right)w|^2dx  \\ &=&
 \int_{\mathbb{R}^n_{k_+}} x_1\cdot \ldots \cdot x_k|\nabla w|^2dx + \frac{1}{4}
 \int_{\mathbb{R}^n_{k_+}} x_1\cdot \ldots \cdot x_k
 \left(\frac{1}{x_1^2}+\ldots+\frac{1}{x_k^2} \right)|w|^2dx  \\
  &+& \frac{1}{2} \int_{\mathbb{R}^n_{k_+}} x_1\cdot \ldots \cdot x_k
 \left( \frac{1}{x_1}, \ldots,
 \frac{1}{x_k} \right)\nabla w^2  dx.
\end{eqnarray*}
By partial integration, we see that the last term is equal to
zero. If the second term is expressed in terms of $u$, we see that
it is equal to the Hardy term
\[
\frac{1}{4}\int_{\mathbb{R}^n_{k_+}}
 \left(\frac{1}{x_1^2}+\ldots+\frac{1}{x_k^2}\right)|u|^2dx.
\]
By Theorem C, the first term may be estimated from below
by the Sobolev term provided that we can prove the following $L^1$
Hardy inequality.
\[
 C\int_{\mathbb{R}^n_{k_+}} (x_1\cdot \ldots \cdot x_k)^{\frac{n-1}{n-2}}(\frac{1}{x_1^2}+\ldots
 +\frac{1}{x_k^2})^{\frac{1}{2}}|v|dx
 \leq \int_{\mathbb{R}^n_{k_+}} (x_1\cdot \ldots \cdot x_k)^{\frac{n-1}{n-2}}|\nabla v|dx.
\]
To do this we work as  in the previous section,  using the
inequality
\[
\left| \int_{\mathbb{R}^n_{k_+}} \mathrm{div} \mathbf{F} |v|dx \right| \leq
\int_{\mathbb{R}^n_{k_+}}|\mathbf{F}|\nabla v|dx,
\]
with the proper choice of vector field,  which turns out to be
\[
  \mathbf{F}=
 (x_1\cdot \ldots \cdot x_k)^{\tau}
\left(\frac{1}{x_1^2}+\ldots+\frac{1}{x_k^2}\right)^{\beta}\left(
 \frac{1}{x_1},\ldots,\frac{1}{x_k}\right),
\]
where
\[
  \tau = \frac{n-1}{n-2} \quad \textrm{ and } \quad \beta = -\frac{1}{2}.
\]
We immediately see that
 $|\mathbf{F}|=\phi^{2\tau}=(x_1\cdot \ldots \cdot x_k)^{\frac{n-1}{n-2}}$.  Also,
\begin{eqnarray*}
 \mathrm{div}  \mathbf{F}&=& -(x_1\cdot \ldots \cdot x_k)^{\tau}
 \left(\frac{1}{x_1^2}+\ldots+\frac{1}{x_k^2} \right)^{\beta +
 1}\\
 &+& \tau (x_1\cdot \ldots \cdot x_k)^{\tau}
 \left(\frac{1}{x_1^2}+\ldots+\frac{1}{x_k^2} \right)^{\beta+1} \\
 &+&
 -2\beta \left(\frac{1}{x_1^4}+\ldots+\frac{1}{x_k^4} \right) \cdot
 (x_1\cdot \ldots \cdot x_k)^{\tau}
 \left(\frac{1}{x_1^2}+\ldots+\frac{1}{x_k^2} \right)^{\beta-1} \\
\end{eqnarray*}
Since $\tau-1>0$ and the last term is positive, we get the result.

\end{proof}

\end{document}